\documentclass{article}
\usepackage{latexsym,amsfonts,palatino,eulervm,color,amsmath,amssymb,amsthm,graphics,graphicx,comment,float,authblk}
\usepackage[hang,flushmargin]{footmisc}

\newtheorem{theorem}{Theorem}[section]
\newtheorem{lemma}[theorem]{Lemma}
\newtheorem{conjecture}[theorem]{Conjecture}
\newtheorem{definition}[theorem]{Definition}

\newtheorem{problem}[theorem]{Problem}
\newtheorem{example}[theorem]{Example}

\def\n{\noindent}

\def\F{\mathbb{F}}

\def\lcm{\mathrm{lcm}}

\def\rank{\mathrm{rank}}

\title{Groups of singular alternating sign matrices}
\author[1]{Cian O'Brien}
\author[2]{Rachel Quinlan}
\affil[1]{Department of Mathematics and Computer Studies, Mary Immaculate College}
\affil[2]{School of Mathematical and Statistical Sciences, University of Galway}

\begin{document}
\maketitle
\begin{abstract}
We investigate multiplicative groups consisting entirely of singular alternating sign matrices (ASMs), and present several constructions of such groups.  It is shown that every finite group is isomorphic to a group of singular ASMs, with a singular idempotent ASM as its identity element. The relationship between the size, the rank, and the possible multiplicative orders of singular ASMs is explored.

\end{abstract}

\n\emph{Keywords}: Alternating sign matrix, matrix group

\n\emph{MSC}: 15A30, 15B36, 05B20

\section{Introduction}

An \emph{alternating sign matrix (ASM)} is a matrix with entries from $\{0, +1, -1\}$ for which the non-zero entries in each row and column alternate in sign, starting and ending with $+1$. Immediate consequences of this definition are that ASMs are square, and each row and column sums to 1. Alternating sign matrices were first investigated by Mills, Robbins, and Rumsey \cite{MRR}, in a context arising from the classical theory of determinants. Connections to fields such as statistical mechanics \cite{Zeilberger} and enumerative combinatorics \cite{Kuperberg} were subsequently discovered, and ASMs continue to attract sustained interest from diverse viewpoints. We refer to Bressoud's book \cite{Bressoud}, for a comprehensive account of the emergence of attention to ASMs and the mathematical developments that ensued. 

A recurrent theme in the study of ASMs is their occurrence, in independent contexts, as natural generalizations of permutation matrices. This invites the question of whether and how familiar themes in the study of permutations can be applied or adapted to ASMs. For example, ASMs first emerged in the definition of the $\lambda$-determinant of a square matrix, which involves adapting the technique of Dodgson condensation by replacing the usual $2\times 2$ determinant with a version involving a parameter $\lambda$. Alternating sign matrices play the role for the $\lambda$-determinant that permutations do for the special case of the classical determinant, which arises if the value of $\lambda$ is set to $-1$. Lascoux and Sch\"{u}tzenberger \cite{bruhat} showed that the set of $n\times n$ ASMs is the unique minimal lattice extension of the set of $n\times n$ permutation matrices under the \emph{Bruhat partial order}. An extension of the concept of Latin squares, which arise by replacing permutation matrices with ASMs, is investigated in \cite{ashl}.

Under matrix multiplication, the permutation matrices of size $n\times n$ comprise a subgroup of the general linear group. This is a property that does not extend at all, within the class of ASMs. It is shown in \cite{finiteASM} that if a pair of nonsingular ASMs are inverses of each other, then they are permutation matrices. It follows that the only subgroups of a general linear group whose elements are all ASMs are groups of permutation matrices. However, it is possible for a set of singular ASMs to form a group under matrix multiplication, with an identity element that is an idempotent ASM. The purpose of this article is to initiate an exploration of this phenomenon. We present examples and constructions, and propose some questions for further investigations. We show in particular that the symmetric group $S_n$, and hence any group of order $n$, can be realized in multiple ways as multiplicative groups of singular ASMs.

\section{Groups of Singular ASMs and Associated Linear Groups} \label{groups}
Throughout this paper, we say that a (non-zero) square matrix $A$ has finite order if there exists an integer $k \geq 1$ such that $A^{k+1} = A$. We say that $A$ has order $k$ if $k$ is the minimal such integer. We let $SA_n$ denote the set of $n \times n$ singular ASMs.

\begin{example}\label{5x5} Consider the following matrices in $SA_5$.
$$
\mathcal{E}_1 = \left(\begin{array}{rrrrr}
0 & 0 & 1 & 0 & 0 \\
0 & 1 & -1 & 0 & 1 \\
0 & 0 & 1 & 0 & 0 \\
1 & 0 & -1 & 1 & 0 \\
0 & 0 & 1 & 0 & 0
\end{array}\right)\hspace{2cm}A= \left(\begin{array}{rrrrr}
0 & 0 & 1 & 0 & 0 \\
1 & 0 & -1 & 1 & 0 \\
0 & 0 & 1 & 0 & 0 \\
0 & 1 & -1 & 0 & 1 \\
0 & 0 & 1 & 0 & 0
\end{array}\right)
$$
$\mathcal{E}_1$ is idempotent (has order 1) since $\mathcal{E}_1^2 = \mathcal{E}_1$. Also $A^2=\mathcal{E}_1$, and $\mathcal{E}_1A=A\mathcal{E}_1=A$, so $A$ has order 2 and $\{\mathcal{E}_1,A\}$ is a group of order 2 under matrix multplication.
\end{example}

We refer to a multiplicative group of singular matrices as a \emph{singular group}, and to a multiplicative group of non-singular matrices as a \emph{linear group}. 
Example \ref{5x5} establishes the existence of a singular group in $SA_5$ that is cyclic of order 2. In this section, we note some general properties of singular groups.

Let $G$ be a singular group of $n\times n$ matrices, with entries in a field $\F$. The identity element $E$ of $G$ is an idempotent matrix, and all elements of $G$ have the same rank $r$, the same rowspace, and the same columnspace. If we write $E=I_n-T$, then it follows from $E^2=E$ that $T^2=T$, and so $T$ is also idempotent. Moreover, for every $A\in G$
$$
EA=A=AE\Longrightarrow (I_n-T)A=A=A(I_n-T)=A\Longrightarrow TA=AT=0.
$$
It follows that $\rank (T)\le n-r$. On the other hand the rank of $T$ is at least $n-r$, since $T+E=I_n$. Thus $T$ has rank $n-r$ and its rowspace and columnspace are exactly the left and right nullspaces of $E$ (and of every element of $G$).

Our next theorem demonstrates a natural correspondence between singular groups and linear groups.

\begin{theorem}\label{gamma}
    Let $G$ be a singular group of matrices in $\mathbb{F}^{n \times n}$ with identity $E = I_n - T$, and let $\Gamma = \{X+T: X\in G\}$. Then $\Gamma$ is a linear group, isomorphic to $G$ via the bijection $X\longleftrightarrow X+T$.
    \end{theorem}

\begin{proof}
    Let $A,B\in G$. Then $(A+T)(B+T) = AB+TB+AT+T^2 = AB+T$. Since $E+T=I_n$, it follows that the mapping from $G$ to $\Gamma$ defined by $X\to X+T$ is a homomorphism of groups, and that its image consists of matrices of rank $n$. It is clearly injective, so its image is isomorphic to $G$.
\end{proof}

We remark that the isomorphism from $\Gamma$ to $G$ of Theorem \ref{gamma} may be described as subtraction of $T$ or as multiplication (on the left or right) by $E$. For every element $Y$ of $\Gamma$ we have $Y-T=YE=EY$.

A possible strategy for constructing a group $G$ of singular ASMs is to start with a particular idempotent $E\in SA_n$ as the identity element of $G$, and try to identify a linear group that is a candidate for $\Gamma$. This would require that every element of $\Gamma$ commutes with $E$, but we may first suppose that $\Phi$ is a linear group with the property that  $E\Phi =\{EH: H\in\Phi\}$ is a set of ASMs. For $E\Phi$ to be a group requires that $EH$ and $E$ have the same row space for every $H\in\Phi$ (that they have the same column space is automatic, since each $H\in\Phi$ is nonsingular). Candidates for $\Phi$ are groups of matrices that preserve the alternating sign structure when multiplying $E$ on the right, and also preserve the rowspace of $E$. These are substantial restrictions, and the obvious candidates are likely to be groups of permutation matrices.

Now suppose that $E\Phi$ \emph{is} a group for which $EH$ has the same row space and the same column space as $E$, for every $H\in\Phi$. It follows that $EHE=EH$ for every $H\in\Phi$, since every row of $EH$ belongs to the left 1-eigenspace of $E$. If $K,L\in\Phi$, then $EKEL=EKL$, and it follows that the mapping $H\rightarrow EH$ is a homomorphism of groups. This is not necessarily injective, since $EH=E$ does not require that $H$ is the identity matrix, only that the columnspace of $H-I$ is contained in that of $T$. In any case, $E\Phi$ is isomorphic to a quotient of $\Phi$. 

The linear group $\Phi$ need not be equal to the $\Gamma$ of Theorem \ref{gamma}, and it is possible that $E\Phi$ is a group but that $\Phi E = \{HE: H\in\Phi\}$ is not. The unique linear group $\Gamma$ determined by the singular group $E\Phi$ is $\{EH+T:H\in\Phi\}$, where $T=I-E$ and $I$ is the identity matrix. In general, the matrices $EH+T$ and $H$ differ by an element of the right nullspace of $E$. 

We now consider the idempotent $\mathcal{E}_1$ of Example \ref{5x5}. Permuting the last two columns of this matrix, or the first two, preserves the alternating sign property. Neither preserves the row space of $\mathcal{E}_1$, which has a basis $\{[0\ 0\ 1\ 0\ 0],\ [0\ 1\ 0\ 0\ 1],\ [1\ 0\ 0\ 1\ 0]\}$. To preserve the alternating sign structure \emph{and} the row space, we can apply both of these column transpositions to obtain the matrix $A$ of Example \ref{5x5}. As a candidate for $\Phi$, this suggests the copy of $C_2$ generated by the permutation
$$
P=\left(\begin{array}{rrrrr}
0 & 1 & 0 & 0 & 0 \\
1 & 0 & 0 & 0 & 0 \\
0 & 0 & 1 & 0 & 0 \\
0 & 0 & 0 & 0 & 1 \\
0 & 0 & 0 & 1 & 0 \\
\end{array}\right).
$$
With $T = I_5-\mathcal{E}_1$, the corresponding linear group $\Gamma$ is the copy of $C_2$ generated by
$$
A'=A+T=\left(\begin{array}{rrrrr}
1 & 0 & 0 & 0 & 0 \\
1 & 0 & 0 & 1 & -1 \\
0 & 0 & 1 & 0 & 0 \\
-1 & 1 & 0 & 0 & 1 \\
0 & 0 & 0 & 0 & 1 \\
\end{array}\right).
$$
It is easily confirmed that $\mathcal{E}_1(A'-P)=0$. 

As an alternative to seeking a linear group $\Phi$ that acts on the columnspace of $\mathcal{E}_1$, one could look for a linear group $\Theta$ for which $\{HE: H\in\Theta\}$ is a singular group of ASMs. This requires that each $HE$ is an ASM with the same rowspace as $E$. Switching rows 2 and 4 of $\mathcal{E}_1$ produces $A$, so a candidate for $\Theta$ is the copy of $C_2$ generated by the transposition 
$$
Q=\left(\begin{array}{rrrrr}
1 & 0 & 0 & 0 & 0 \\
0 & 0 & 0 & 1 & 0 \\
0 & 0 & 1 & 0 & 0 \\
0 & 1 & 0 & 0 & 0 \\
0 & 0 & 0 & 0 & 1 \\
\end{array}\right).
$$
We note that $(Q-A')\mathcal{E}_1=0$. Neither $\mathcal{E}_1Q$ nor $P\mathcal{E}_1$ is an ASM, although each generates a singular group of order 2. The linear group $\Gamma$ of order 2 is generated by $A'$, which satisfies
$$
Q\mathcal{E}_1=A'\mathcal{E}_1=\mathcal{E}_1A'=\mathcal{E}_1P.
$$

\section{Idempotent ASMs}\label{section idempotent}
The loose strategy outlined in Section \ref{groups} for constructing singular groups of ASMs relies on being able to identify idempotent ASMs to act as identity elements in such groups. In this section, we take some initial steps in this direction by describing all idempotent ASMs of nullity at most 2.

\begin{definition}
  An alternating sign matrix $A$ is said to be \emph{in reduced form} if for every $i$ with $A_{ii}=1$, either row $i$ or column $i$ contains a negative entry.
\end{definition}

Every ASM $A$ has a unique submatrix that is an ASM in reduced form, called \emph{the} reduced form of $A$, which can be obtained from $A$ by repeating the step of deleting row $i$ and column $i$ for an index $i$ with the property that a 1 in position $(i,i)$ is the only non-zero entry in both row $i$ and column $i$. We say that $A$ is a \emph{diagonal extension} of its reduced form. Deleting such a diagonal entry along with its row and column preserves the property of being an ASM, the nullity, and the minimum polynomial. In particular, if $A$ is an idempotent ASM, then so is its reduced form. It follows that in order to describe all idempotent ASMs with specified nullity, it is enough to identify those that are in reduced form.

The only $n\times n$ idempotent matrix with nullity 0 is the identity matrix, which is an ASM. The only groups of ASMs of nullity 0 are groups of permutation matrices {\cite{finiteASM}}. The next lemma shows that there are no idempotent ASMs of nullity 1. 

\begin{lemma}\label{rank n-2}   
    Let $E \in SA_n$ be idempotent. Then the nullity of $E$ is at least 2. 
\end{lemma}

\begin{proof}
  It is sufficient to prove the lemma in the case where $E$ is in reduced form, and we make this assumption. It follows immediately that $E_{11}=E_{nn}=0$, since the first row, last row, first column and last column of an ASM each contain only one non-zero entry.
  Moreover $E_{1n}=0$ also, since $E_{1n}=(E^2)_{1n}=1$ would require $E_{nn}=1$. Similarly $E_{n1}=0$, so the four corner entries of $E$ are zero. Let $j_1$ be the position of the 1 in row 1 of $E$. Then $2\le j_1\le n-1$. Since $E^2=E$, row $j_1$ of $E$ is a repeat of row 1. Similarly, if $j_2$ is the position of the 1 in row $n$, then row $j_2$ and row $n$ are duplicates. This establishes the existence of two independent elements of the left nullspace of $E$, hence the nullity of $E$ is at least 2.
\end{proof}

As a consequence of Lemma \ref{rank n-2}, we note that no ASM of finite multiplicative order has exactly one negative entry. If $A$ is a $n\times n$ ASM with exactly one negative entry, then the row and column of the negative entry each have exactly three non-zero entries, and every other row and column has a 1 as its unique non-zero entry. There exist permutation matrices $P$ and $Q$ for which $PAQ$ is a block diagonal matrix with the unique non-permutation $3\times 3$ ASM as its upper left block, and the $(n-3)\times (n-3)$ identity matrix as its lower right block. It follows that $A$ has rank $n-1$, and hence that $A$ cannot belong to a singular group of ASMs, since no idempotent ASM has rank $n-1$.

We now proceed to the main result of this section, that all idempotent ASMs of nullity 2 are diagonal extensions of the matrix $\mathcal{E}_1$ of Example \ref{5x5}, or of its transpose.

 \begin{theorem}\label{nullity2} Let $E$ be an idempotent ASM of nullity 2 that is in reduced form. Then $E$ is either the $5\times 5$ matrix $\mathcal{E}_1$ of Example \ref{5x5}, or its transpose.
 \end{theorem}
 
\begin{proof}
Let $j_1,j_2,i_1,i_2$ respectively be the positions of the 1 in row 1, row $n$, column 1, and column $n$ of $E$. These indices are all in the range 2 to $n-1$, since the corner entries are all $0$, as noted in the proof of Lemma \ref{rank n-2}. Looking at the corner entries of $E^2$, we also note that $\{i_1,i_2\}\cap\{j_1,j_2\}=\emptyset$.
In particular $i_1\neq j_1$, and we assume that $i_1<j_1$, replacing $E$ with its transpose if necessary.

As in Section \ref{groups}, we write $T$ for the matrix $I_n-E$, which has rank 2. Row 1 of $T$ has 1 in position 1, $-1$ in position $j_1$, and 0 in all other positions. Row $n$ has $-1$ in position $j_2$, 1 in position $n$, and 0 in all other positions. Since row 1 and row $n$ are linearly independent, they span the rowspace of $T$. Column 1 and column $n$ of $T$ each has only two nonzero entries, in positions $1$ and $i_1$ in column 1, and in positions $i_2$ and $n$ in column $n$. Every linear combination of row 1 and row $n$ has a nonzero entry either in position 1 or in position $n$. It follows that the only nonzero rows of $T$ are rows 1, $n$, $i_1$ and $i_2$. Similarly, the only non-zero columns of $T$ are columns $1$, $n$, $j_1$ and $j_2$. The condition that $E$ is in reduced form means that either row $i$ or column $i$ of $T$ includes a nonzero entry, for $i$ in the range $1$ to $n$. It follows that $n\le 6$, and $n=6$ if and only if $i_1\neq i_2$ and $j_1\neq j_2$. 

Suppose that $i_1\neq i_2$. Then row $i_1$ of $T$ has $0$ as its last entry, so it is a scalar multiple of row 1, with $-1$ as its first entry. The non-zero entries of row $i_1$ are $-1$ and $1$ in positions $1$ and $j_1$ respectively. Since $i_1<j_1$ it follows that the three non-zero entries in row $i_1$ of $E=I_n-T$ occur in the sequence $1,1,-1$, contrary to the hypothesis that $E$ is a ASM. We conclude that $i_1=i_2$, and that $j_2<i_1<j_1$. It follows that $n=5$ and that $(1,2,3,4,5)=(1,j_2,i_1,j_1,5)$. We conclude that 
$$
E= \left(\begin{array}{rrrrr}
            0 & 0 & 0 & 1 & 0 \\
            0 & 1 & 0 & 0 & 0 \\
            1 & -1 & 1 & -1 & 1 \\
            0 & 0 & 0 & 1 & 0 \\
            0 & 1 & 0 & 0 & 0 \\
\end{array}\right).
$$
Thus $E$ is the transpose of $\mathcal{E}_1$. The alternative conclusion, $E=\mathcal{E}_1$, would result by the same argument from the assumption that $i_1>j_1$.
\end{proof}

\section{Groups in $SA_n$ for $n \leq 7$}\label{n leq 7}
As previously discussed, $\mathcal{E}_1$ and its transpose are the only idempotent matrices in $SA_n$ for $n \leq 5$. By exhaustive search for $n \leq 5$, the only other matrices of finite order in $SA_n$ are $A$ from Example \ref{5x5} and its transpose. So for $n \leq 5$, the only multiplicative groups in $SA_n$ are $\{\mathcal{E}_1 \}, \{\mathcal{E}_1 ,A\}, \{\mathcal{E}_1 ^t\},$ and $\{\mathcal{E}_1 ^t,A^t\}$.

For $n=6$, we have seen that each idempotent matrix of nullity 2 in $SA_6$ has $\mathcal{E}_1$ or $\mathcal{E}_1^t$ as its reduced form. If we wish to construct a $6\times6$ idempotent ASM from $\mathcal{E}_1$ or $\mathcal{E}_1^t$, there are 6 choices for where to insert the extra row and column with $+1$ on the diagonal as the only non-zero entry. This means there are 12 idempotent matrices of nullity 2 in $SA_6$. Exhaustive search confirms that these are the only idempotents in $SA_6$, and that there are 28 other matrices of finite order in $SA_6$, each of which squares to one of these idempotent matrices. Each of the 12 idempotent matrices is the square of at least one of these order 2 ASMs, with some being the square of two or four of them.

\begin{example}\label{6x6} Consider the following $6 \times 6$ singular ASMs of rank 4.
$$
E = \left(\begin{array}{rrrrrr}
0 & 0 & 0 & 1 & 0 & 0 \\
0 & 1 & 0 &-1 & 0 & 1 \\
0 & 0 & 1 & 0 & 0 & 0 \\
0 & 0 & 0 & 1 & 0 & 0 \\
1 & 0 & 0 &-1 & 1 & 0 \\
0 & 0 & 0 & 1 & 0 & 0 
\end{array}\right)\hspace{0.1cm}
A = \left(\begin{array}{rrrrrr}
0 & 0 & 0 & 1 & 0 & 0 \\
1 & 0 & 0 &-1 & 1 & 0 \\
0 & 0 & 1 & 0 & 0 & 0 \\
0 & 0 & 0 & 1 & 0 & 0 \\
0 & 1 & 0 &-1 & 0 & 1 \\
0 & 0 & 0 & 1 & 0 & 0 
\end{array}\right)
\hspace{0.1cm}B = \left(\begin{array}{rrrrrr}
0 & 0 & 1 & 0 & 0 & 0 \\
0 & 1 &-1 & 0 & 0 & 1 \\
0 & 0 & 0 & 1 & 0 & 0 \\
0 & 0 & 1 & 0 & 0 & 0 \\
1 & 0 &-1 & 0 & 1 & 0 \\
0 & 0 & 1 & 0 & 0 & 0 
\end{array}\right)
$$
$$
C = \left(\begin{array}{rrrrrr}
0 & 0 & 1 & 0 & 0 & 0 \\
1 & 0 &-1 & 0 & 1 & 0 \\
0 & 0 & 0 & 1 & 0 & 0 \\
0 & 0 & 1 & 0 & 0 & 0 \\
0 & 1 &-1 & 0 & 0 & 1 \\
0 & 0 & 1 & 0 & 0 & 0 
\end{array}\right)\hspace{0.1cm}D = \left(\begin{array}{rrrrrr}
0 & 0 & 0 & 1 & 0 & 0 \\
0 & 0 & 1 & 0 & 0 & 0 \\
0 & 1 & 0 &-1 & 0 & 1 \\
0 & 0 & 0 & 1 & 0 & 0 \\
1 & 0 & 0 &-1 & 1 & 0 \\
0 & 0 & 0 & 1 & 0 & 0 
\end{array}\right)
$$
The matrix $E$ is idempotent and has the same rowspace and columnspace as $A,B,C$ and $D$, so $EX=XE=X$ for $X\in\{A,B,C,D\}$. Furthermore $A^2 = B^2 = C^2 = D^2 = E$, so each of $A,B,C,D$ generates a copy of $C_2$ in $SA_6$. The matrices $A$ and $B$ satisfy $AB=BA=C$ and hence generate a copy of $C_2\times C_2$ in $SA_6$. The group generated by $D$ is maximal as a multiplicative group in $SA_6$.
\end{example}

Note that $A$ has a $5\times5$ ASM of order 2 as its reduced form (resulting from the deletion of row and column 3), whereas $B$, $C$, and $D$ are in reduced form. All idempotent elements of $SA_6$ are diagonal extensions of those in $SA_5$. However, Example \ref{6x6} shows that a singular group of ASMs may contain some elements that are in reduced form, and some that are not.

Let $G$ be the Klein 4-group generated by $A$ and $B$. To consider how $G$ might be constructed from $E$ according to the scheme outlined in Section \ref{groups}, we note that $A$ can be obtained from $E$ by switching columns 1 and 6 and switching columns 2 and 5, and that $B$ can be obtained from $E$ by switching columns 3 and 4. These adjustments from $E$ to $A$ and $B$ can be achieved by multiplying $E$ on the right by permutation matrices corresponding to the elements $(1 \ 6)(2 \ 5)$ and $(3 \ 4)$ of $S_6$, which generate a linear group $\Phi$ isomorphic to $C_2\times C_2$. The elements of $\Phi$ do not commute with $E$. The linear group $\Gamma$ corresponding to $G$ includes the following non-permutation matrices, obtained as in Theorem \ref{gamma} by adding $T=I-E$ to $A,B,C$ respectively.

$$
\left(\begin{array}{rrrrrr}
1 & 0 & 0 & 0 & 0 & 0 \\
1 & 0 & 0 & 0 & 1 & -1 \\
0 & 0 & 1 & 0 & 0 & 0 \\
0 & 0 & 0 & 1 & 0 & 0 \\
-1 & 1 & 0 & 0 & 0 & 1 \\
0 & 0 & 0 & 0 & 0 & 1
\end{array}\right), 
\left(\begin{array}{rrrrrr}
1 & 0 & 1 & -1 & 0 & 0 \\
0 & 1 & -1 & 1 & 0 & 0 \\
0 & 0 & 0 & 1 & 0 & 0 \\
0 & 0 & 1 & 0 & 0 & 0 \\
0 & 0 & -1 & 1 & 1 & 0 \\
0 & 0 & 1 & -1 & 0 & 1
\end{array}\right),
\left(\begin{array}{rrrrrr}
1 & 0 & 1 & -1 & 0 & 0 \\
1 & 0 & -1 & 1 & 1 & -1 \\
0 & 0 & 0 & 1 & 0 & 0 \\
0 & 0 & 1 & 0 & 0 & 0 \\
-1 & 1 & -1 & 1 & 0 & 1 \\
0 & 0 & 1 & -1 & 0 & 1

\end{array}\right)
$$
As in Example \ref{5x5}, $\Gamma$ is not a group of permutation matrices, but there is a group $\Phi$ of permutation matrices for which $G=E\Phi = E\Gamma$. Unlike Example \ref{5x5}, a linear group $\Theta$ for which $\Theta E = G$ cannot consist of permutation matrices, since $B$ cannot be obtained from $E$ by permutation of the rows. An example of such a group $\Theta$ (different from $\Gamma$ but also isomorphic to $C_2\times C_2$) is generated by the matrix of the transposition $(2 \ 5)$ and the matrix $P$ below, with $B=PE$.
$$
P=\left(\begin{array}{rrrrrr}
0 & 0 & 1 & 0 & 0 & 0 \\
1 & 1 & -1 & 0 & 0 & 0 \\
1 & 0 & 0 & 0 & 0 & 0 \\
-1 & 0 & 1 & 1 & 0 & 0 \\
1 & 0 & -1 & 0 & 1 & 0 \\
0 & 0 & 1 & 1 & 0 & -1
\end{array}\right).
$$

For $n = 7$, exhaustive search yields singular ASMs of order 1, 2, 3, 4, and 6, as summarised in the following table.
\[\begin{array}{c|c}
\text{Order}		& \text{Number of Singular $7\times7$ ASMs}\\
\hline
1			& 42	\\
2			& 270 \\
3			& 36	\\
4			& 32	\\
6			& 12	\\
\end{array}\]

If we wish to construct a $7\times7$ idempotent ASM from $\mathcal{E}_1$ or $\mathcal{E}_1^t$, there are ${7 \choose 2} = 21$ choices for where to insert the extra rows and columns with $+1$ on the diagonal and zero elsewhere. This means that all 42 idempotent matrices in $SA_7$ have $\mathcal{E}_1$ or $\mathcal{E}_1^t$ as their reduced forms and all have nullity 2.

Using ASMs in $SA_7$, we can form cyclic groups of order 1, 2, 3, 4, and 6, as well as non-Abelian groups of order 6, 8, and 12.

\begin{example}\label{2 generators} $A$ has order 3 and $B$ has order 2, with $E$ acting as an identity for both.
\[E = \left(\begin{array}{rrrrrrr}
0 & 0 & 1 & 0 & 0 & 0 & 0 \\
0 & 1 &-1 & 0 & 0 & 0 & 1 \\
0 & 0 & 1 & 0 & 0 & 0 & 0 \\
0 & 0 & 0 & 1 & 0 & 0 & 0 \\
0 & 0 & 0 & 0 & 1 & 0 & 0 \\
1 & 0 &-1 & 0 & 0 & 1 & 0 \\
0 & 0 & 1 & 0 & 0 & 0 & 0 
\end{array}\right)\hspace{0.1cm}A = \left(\begin{array}{rrrrrrr}
0 & 0 & 0 & 0 & 1 & 0 & 0 \\
0 & 1 & 0 & 0 &-1 & 0 & 1 \\
0 & 0 & 0 & 0 & 1 & 0 & 0 \\
0 & 0 & 1 & 0 & 0 & 0 & 0 \\
0 & 0 & 0 & 1 & 0 & 0 & 0 \\
1 & 0 & 0 & 0 &-1 & 1 & 0 \\
0 & 0 & 0 & 0 & 1 & 0 & 0 
\end{array}\right)\hspace{0.1cm}
B = \left(\begin{array}{rrrrrrr}
0 & 0 & 0 & 0 & 1 & 0 & 0 \\
0 & 1 & 0 & 0 &-1 & 0 & 1 \\
0 & 0 & 0 & 0 & 1 & 0 & 0 \\
0 & 0 & 0 & 1 & 0 & 0 & 0 \\
0 & 0 & 1 & 0 & 0 & 0 & 0 \\
1 & 0 & 0 & 0 &-1 & 1 & 0 \\
0 & 0 & 0 & 0 & 1 & 0 & 0 
\end{array}\right)\]
$\langle A, B \rangle$ has order 6, and $AB \not= BA$.
\end{example}

\begin{example}\label{order 4} $A$ has order 4 and $B$ has order 2, with $E$ acting as an identity for both.
$$E = \left(\begin{array}{rrrrrrr}
0 & 0 & 0 & 1 & 0 & 0 & 0 \\
0 & 1 & 0 &-1 & 0 & 0 & 1 \\
0 & 0 & 1 & 0 & 0 & 0 & 0 \\
0 & 0 & 0 & 1 & 0 & 0 & 0 \\
1 & 0 & 0 &-1 & 1 & 0 & 0 \\
0 & 0 & 0 & 0 & 0 & 1 & 0 \\
0 & 0 & 0 & 1 & 0 & 0 & 0 
\end{array}\right)\hspace{0.1cm}
A = \left(\begin{array}{rrrrrrr}
0 & 0 & 0 & 1 & 0 & 0 & 0 \\
0 & 0 & 0 & 0 & 0 & 1 & 0 \\
0 & 1 & 0 &-1 & 0 & 0 & 1 \\
0 & 0 & 0 & 1 & 0 & 0 & 0 \\
0 & 0 & 1 & 0 & 0 & 0 & 0 \\
1 & 0 & 0 &-1 & 1 & 0 & 0 \\
0 & 0 & 0 & 1 & 0 & 0 & 0 
\end{array}\right)\hspace{0.1cm}
B = \left(\begin{array}{rrrrrrr}
0 & 0 & 0 & 1 & 0 & 0 & 0 \\
1 & 0 & 0 &-1 & 1 & 0 & 0 \\
0 & 0 & 1 & 0 & 0 & 0 & 0 \\
0 & 0 & 0 & 1 & 0 & 0 & 0 \\
0 & 1 & 0 &-1 & 0 & 0 & 1 \\
0 & 0 & 0 & 0 & 0 & 1 & 0 \\
0 & 0 & 0 & 1 & 0 & 0 & 0 
\end{array}\right)$$
$\langle A, B \rangle$ has order 8, and $AB \not= BA$.
\end{example}

\begin{example} $A$ has order 6 and $B$ has order 2, with $E$ acting as an identity for both.
$$E = \left(\begin{array}{rrrrrrr}
0 & 0 & 0 & 0 & 1 & 0 & 0 \\
0 & 1 & 0 & 0 &-1 & 0 & 1 \\
0 & 0 & 1 & 0 & 0 & 0 & 0 \\
0 & 0 & 0 & 1 & 0 & 0 & 0 \\
0 & 0 & 0 & 0 & 1 & 0 & 0 \\
1 & 0 & 0 & 0 &-1 & 1 & 0 \\
0 & 0 & 0 & 0 & 1 & 0 & 0 
\end{array}\right)\hspace{0.1cm}
A = \left(\begin{array}{rrrrrrr}
0 & 0 & 0 & 1 & 0 & 0 & 0 \\
1 & 0 & 0 &-1 & 0 & 1 & 0 \\
0 & 0 & 0 & 0 & 1 & 0 & 0 \\
0 & 0 & 1 & 0 & 0 & 0 & 0 \\
0 & 0 & 0 & 1 & 0 & 0 & 0 \\
0 & 1 & 0 &-1 & 0 & 0 & 1 \\
0 & 0 & 0 & 1 & 0 & 0 & 0 
\end{array}\right)\hspace{0.1cm}
B = \left(\begin{array}{rrrrrrr}
0 & 0 & 0 & 0 & 1 & 0 & 0 \\
1 & 0 & 0 & 0 &-1 & 1 & 0 \\
0 & 0 & 0 & 1 & 0 & 0 & 0 \\
0 & 0 & 1 & 0 & 0 & 0 & 0 \\
0 & 0 & 0 & 0 & 1 & 0 & 0 \\
0 & 1 & 0 & 0 &-1 & 0 & 1 \\
0 & 0 & 0 & 0 & 1 & 0 & 0 
\end{array}\right)$$
$\langle A, B \rangle$ has order 12, and $AB \not= BA$.
\end{example}

Exhaustive search shows that there is no group of greater order than 12 in $SA_7$. Note in particular that the group in Example \ref{2 generators} is isomorphic to $S_3$, but there is no group isomorphic to $S_n$ in $SA_7$ for any $n \geq 3$.

All examples of finite order ASMs in $SA_n$ discussed so far have had only two negative entries, both occurring in the same row or column. For $n \geq 7$, we also find finite order ASMs in $SA_n$ containing more than two negative entries. All such examples for $n = 7$ contain three negative entries, and square to an idempotent matrix. No such example can be multiplied by any other finite order ASM in $SA_7$, besides itself or its square, to give another finite order singular ASM.

\begin{example}\label{3 negatives}
$A$ has order 2, and $E = A^2$ acts as its identity.
$$A = \left(\begin{array}{rrrrrrr}
0 & 0 & 0 & 0 & 1 & 0 & 0 \\
0 & 0 & 1 & 0 &-1 & 1 & 0 \\
0 & 0 & 0 & 0 & 1 & 0 & 0 \\
1 & 0 &-1 & 1 & 0 & 0 & 0 \\
0 & 0 & 1 & 0 & 0 & 0 & 0 \\
0 & 1 & 0 & 0 &-1 & 0 & 1 \\
0 & 0 & 0 & 0 & 1 & 0 & 0 
\end{array}\right)\hspace{1cm}E = \left(\begin{array}{rrrrrrr}
0 & 0 & 1 & 0 & 0 & 0 & 0 \\
0 & 1 &-1 & 0 & 0 & 0 & 1 \\
0 & 0 & 1 & 0 & 0 & 0 & 0 \\
1 & 0 &-1 & 1 & 0 & 0 & 0 \\
0 & 0 & 0 & 0 & 1 & 0 & 0 \\
0 & 0 & 0 & 0 & 0 & 1 & 0 \\
0 & 0 & 1 & 0 & 0 & 0 & 0 
\end{array}\right)
$$
\end{example}

\section{Groups in $SA_n$ of Maximal Rank}\label{max rank}

In Section \ref{section idempotent}, it was shown that the maximum possible rank of an idempotent in $SA_n$ is $n-2$, and that every example achieving this bound is a diagonal extension of $\mathcal{E}_1$ (of Example \ref{5x5}) or its transpose. In this section we outline a general construction of singular groups of ASMs of rank $n-2$ and size $n\times n$, for $n\ge 5$.

\begin{theorem}\label{block perm orders}
Let $P$ be a $n \times n$ permutation matrix of multiplicative order $k$. Then the following $(n+4)\times (n+4)$ matrices $A$ and $B$ are singular ASMs of finite order. The order of $A$ is $k$, and the order of $B$ is $2k$ if $k$ is odd, or $k$ if $k$ is even.

\[A = \left(\begin{array}{cc|ccccr|cc}
0&0  &  0&0&\dots&0&1  &  0&0\\
0&1  &  0&0&\dots&0&-1  &  0&1\\
\hline
0&0  &    &  & &  &    &  0&0\\
\vdots&\vdots &    &  &P& &    &  \vdots&\vdots\\
0&0  &    &  & &  &    &  0&0\\
\hline
1&0  &  0&0&\dots&0&-1 &  1&0\\
0&0  &  0&0&\dots&0&1  &  0&0
\end{array}\right)
\hspace{0.5cm}
B= \left(\begin{array}{cc|ccccr|cc}
0&0  &  0&0&\dots&0&1  &  0&0\\
1&0  &  0&0&\dots&0&-1  &  1&0\\
\hline
0&0  &    &  & &  &    &  0&0\\
\vdots&\vdots &    &  &P& &    &  \vdots&\vdots\\
0&0  &    &  & &  &    &  0&0\\
\hline
0&1  &  0&0&\dots&0&-1 &  0&1\\
0&0  &  0&0&\dots&0&1  &  0&0
\end{array}\right)\]
\end{theorem}

\begin{proof}
Since $P$ has exactly one 1 in its final column, with all other entries 0, the corresponding columns of $A$ and $B$ satisfy the necessary conditions to be ASMs. Since the first and last rows of $A$ and $B$ are equal, the matrices are singular.

Let $E_{ij}$ be the $2\times2$ matrix with a 1 in the $(i,j)$-position and 0s elsewhere, and let $X_i$ be the $2\times n$ matrix with 1 in the $(1,i)$-position, $-1$ in the $(2,i)$-position, and 0 elsewhere. Then $A$ and $B$ can be written as follows.
\[A = \left(\begin{array}{c|c|c}
E_{22}&X_n&E_{22}\\
\hline
0&P&0\\
\hline
E_{11}&-X_n&E_{11}
\end{array}\right)
\hspace{0.5cm}
B= \left(\begin{array}{c|c|c}
E_{21}&X_n&E_{21}\\
\hline
0&P&0\\
\hline
E_{12}&-X_n&E_{12}
\end{array}\right)\]

Since $E_{ab}E_{cd}$ is $E_{ad}$ if $b=c$ and 0 otherwise, we have the following. Let $\mu(k) = \frac{1}{2}(1+(-1)^k)$.
\[A^2 = \left(\begin{array}{c|c|c}
E_{22}&X_nP&E_{22}\\
\hline
0&P^2&0\\
\hline
E_{11}&-X_nP&E_{11}
\end{array}\right)
,\dots,\hspace{0.1cm}
A^k= \left(\begin{array}{c|c|c}
E_{22}&X_nP^{k-1}&E_{22}\\
\hline
0&I_n&0\\
\hline
E_{11}&-X_nP^{k-1}&E_{11}
\end{array}\right)
,\hspace{0.1cm}
A^{k+1}= \left(\begin{array}{c|c|c}
E_{22}&X_nI_n&E_{22}\\
\hline
0&P&0\\
\hline
E_{11}&-X_nI_n&E_{11}
\end{array}\right)\]

\[B^2 = \left(\begin{array}{c|c|c}
E_{22}&X_nP&E_{22}\\
\hline
0&P^2&0\\
\hline
E_{11}&-X_nP&E_{11}
\end{array}\right)
,\dots,\hspace{0.1cm}
B^{k+1}= \left(\begin{array}{c|c|c}
E_{2,1+\mu(k)}&X_n&E_{2,1+\mu(k)}\\
\hline
0&P&0\\
\hline
E_{1,2-\mu(k)}&-X_n&E_{1,2-\mu(k)}
\end{array}\right)
\]

The central block of each power $A^2, A^3, \dots, A^k$ of $A$ is a permutation matrix not equal to $P$, so none of these powers of $A$ are equal to $A$. Since $A^{k+1}=A$, the order of $A$ is $k$ and $A^k$ is the identity element for the group generated by $A$.

Similarly, the least powers of $B$ for which the central block is equal to $P$ are $B^{k+1}$ and $B^{2k+1}$. If $k$ is even, then $\mu(k) = 0$, and $B^{k+1} = B$. If $k$ is odd, then $\mu(k) = 1$, and $B^{k+1} \not= B$. In either case, $B^{2k+1} = B$. So if $k$ is even, the order of $B$ is $k$ and $B^k$ is the identity element for the group generated by $B$. If $k$ is odd, then the order of $B$ is $2k$ and $B^{2k}$ is the identity for this group.
\end{proof}

This construction generates many of the $7 \times 7$ examples outlined at the beginning of this paper, specifically those of order 1, 2, 3, and 6.

\begin{example} The following $7 \times 7$ singular ASMs $A$ and $B$ have order 3 and 6, respectively.
\[A = \left(\begin{array}{cc|ccr|cc}
0&0  &  0&0&1  &  0&0\\
0&1  &  0&0&-1  &  0&1\\
\hline
0&0  &    0& 0 &1    &  0&0\\
0&0 &    1& 0&   0 &  0&0\\
0&0  &    0& 1 & 0&   0&0\\
\hline
1&0  &  0&0&-1 &  1&0\\
0&0  &  0&0&1  &  0&0
\end{array}\right)
\hspace{0.5cm}
B= \left(\begin{array}{cc|ccr|cc}
0&0  &  0&0&1  &  0&0\\
1&0  &  0&0&-1  &  1&0\\
\hline
0&0  &    0& 0 &1    &  0&0\\
0&0 &    1& 0&   0 &  0&0\\
0&0  &    0& 1 & 0&   0&0\\
\hline
0&1  &  0&0&-1 &  0&1\\
0&0  &  0&0&1  &  0&0
\end{array}\right)\]
\end{example}

As in the proof of Theorem \ref{block perm orders}, $E_{ij}$ is the $2 \times 2$ matrix with a 1 in the $(i,j)$-position and 0s elsewhere, and $X_k$ is the $2 \times n$ matrix with 1 in the $(1,k)$-position, $-1$ in the $(2,k)$-position, and 0s elsewhere.

\begin{lemma}\label{same identity}
Let $A$ and $B$ be the following singular ASMs of finite multiplicative order, where $P$ and $Q$ are $n \times n$ permutation matrices, and $\{(E_1,E_2), (F_1,F_2)\} \subseteq \{(E_{11},E_{22}), (E_{12},E_{21})\}$.
\[A = \left(\begin{array}{c|c|c}
E_{2}&X_i&E_{2}\\
\hline
0&P&0\\
\hline
E_{1}&-X_i&E_{1}
\end{array}\right)
\hspace{0.5cm}
B= \left(\begin{array}{c|c|c}
F_{2}&X_j&F_{2}\\
\hline
0&Q&0\\
\hline
F_{1}&-X_j&F_{1}
\end{array}\right)\]
\begin{itemize}
    \item[(a)] The cyclic groups generated by $A$ and $B$ have the same identity element if and only if the non-zero entries in column $i$ of $P$ and column $j$ of $Q$ are in the same position.
    \item[(b)] Suppose the non-zero entries in column $i$ of $P$ and column $j$ of $Q$ are in the same position. Then $AB$ is given by the following, where $(G_1,G_2) \in \{(E_{11},E_{22}), (E_{12},E_{21})\}$ and the non-zero entry in column $k$ of $PQ$ is in the same position as in column $i$ of $P$.
\[AB = \left(\begin{array}{c|c|c}
G_2&X_k&G_2\\
\hline
0&PQ&0\\
\hline
G_1&-X_k&G_1
\end{array}\right)\]
\end{itemize}
\end{lemma}

\begin{proof}\phantom{hello\\}
\begin{itemize}
\item[(a)] Let $a$ and $b$ be the multiplicative orders of $A$ and $B$, respectively. From Theorem \ref{block perm orders}, the identities for $\langle A \rangle$ and $\langle B \rangle$ are the following.

\[id_A = A^a = \left(\begin{array}{c|c|c}
E_{22}&X_iP^{a-1}&E_{22}\\
\hline
0&I_n&0\\
\hline
E_{11}&-X_iP^{a-1}&E_{11}
\end{array}\right)
\hspace{0.5cm}
id_B = B^b = \left(\begin{array}{c|c|c}
E_{22}&X_jQ^{b-1}&E_{22}\\
\hline
0&I_n&0\\
\hline
E_{11}&-X_jQ^{b-1}&E_{11}
\end{array}\right)\]

So $id_A = id_B$ if and only if $X_iP^{a-1} = X_jQ^{b-1}$. Since $P$ is a permutation and $X_i$ has non-zeros only in column $i$, it follows that $X_iP$ has non-zeros only in the column corresponding to the position of the non-zero entry in row $i$ of $P$. Similarly, $X_jQ$ has non-zeros only in the column corresponding to the position of the non-zero entry in row $j$ of $Q$. Multiplying $X_i$ (or $X_j$) by subsequent powers of $P$ (or $Q$) continues to move the only column of non-zero entries according to the underlying permutation.

Since $P^a = I_n$, it follows that $P^{a-1}$ is the inverse of $P$. Therefore the position of the non-zero entry in row $i$ of $P^{a-1}$ corresponds to the position of the non-zero entry in column $i$ of $P$. Similarly, the position of the non-zero entry in row $j$ of $Q^{b-1}$ corresponds to the position of the non-zero entry in column $j$ of $Q$. So $X_iP^{a-1} = X_jQ^{b-1}$ if and only if the non-zero entries in column $i$ of $P$ and column $j$ of $Q$ are in the same position.

\item[(b)] 

\[AB = \left(\begin{array}{c|c|c}
E_{2}F_{2}+E_2F_1&X_iQ&E_{2}F_{2}+E_2F_1\\
\hline
0&PQ&0\\
\hline
E_{1}F_{2}+E_1F_1&-X_iQ&E_{1}F_{2}+E_1F_1
\end{array}\right)\]

By considering each possible combination of $E = (E_1,E_2)$ and $F=(F_1,F_2)$, we see that $(E_{1}F_{2}+E_1F_1, E_{2}F_{2}+E_2F_1)$ is $(E_{11}, E_{22})$ if $E=F$, and is $(E_{12}, E_{21})$ if $E\not=F$. Therefore $E_{1}F_{2}+E_1F_1 = G_1$ and $E_{2}F_{2}+E_2F_1 = G_2$. After right multiplication by $Q$, the non-zero entry in column $i$ of $P$ is sent to the same column as $X_i$. Therefore $X_iQ = X_k$, and $AB$ is as stated.
\end{itemize}
\end{proof}

\begin{theorem}\label{isomorphic n+4}
$S_n$ is isomorphic to a group of singular ASMs in $SA_{n+4}$.
\end{theorem}
\begin{proof}
We represent the $n$-cycle $(1,2,\dots,n)$ and the transposition $(1,n)$ by the following $(n+4)\times(n+4)$ singular ASMs $S$ and $T$, respectively.

\[S = \left(\begin{array}{cc|ccccc|cc}
0&0  &  0&0&\dots&0&1  &  0&0\\
0&1  &  0&0&\dots&0&-1  &  0&1\\
\hline
0&0  &   		&  &  &  &1   	&  0&0\\
0&0  &   		1&  &  &  &    	&  0&0\\
\vdots&\vdots&    	&  1&  &  &    	&  \vdots&\vdots\\
0&0  &    		&  &  \ddots&  &    	&  0&0\\
0&0  &    		&  &  &  1&    	&  0&0\\
\hline
1&0  &  0&0&\dots&0&-1 &  1&0\\
0&0  &  0&0&\dots&0&1  &  0&0
\end{array}\right)
\hspace{0.5cm}
T = \left(\begin{array}{cc|ccccc|cc}
0&0  &  0&0&\dots&0&1  &  0&0\\
1&0  &  0&0&\dots&0&-1  &  1&0\\
\hline
0&0  	&			   &  &  &  & 1   		&  0&0\\
0&0  	&			   & 1 &  &  &    		&  0&0\\
\vdots&\vdots &  		   &  & \ddots &  &    		&  \vdots&\vdots\\
0&0  	&  			   &  &  & 1 &    		&  0&0\\
0&0  	&  			1 &  &  &  &    		&  0&0\\
\hline
0&1  &  0&0&\dots&0&-1 &  0&1\\
0&0  &  0&0&\dots&0&1  &  0&0
\end{array}\right)\]

More concisely, $S$ and $T$ can be represented by the block matrices.

\[S = \left(\begin{array}{c|c|c}
E_{22}&X_n&E_{22}\\
\hline
0&P&0\\
\hline
E_{11}&-X_n&E_{11}
\end{array}\right)
\hspace{0.5cm}
T = \left(\begin{array}{c|c|c}
E_{21}&X_n&E_{21}\\
\hline
0&Q&0\\
\hline
E_{12}&-X_n&E_{12}
\end{array}\right)\]

Since $P$ has order $n$ and $Q$ has order 2, Theorem \ref{block perm orders} implies that $S$ has order $n$ and $T$ has order 2. Since the non-zero entry in column $n$ of $P$ and $Q$ are in the same position, Lemma \ref{same identity} (a) implies that $\langle S \rangle$ and $\langle T \rangle$ have the same identity element. Lemma \ref{same identity} (b) implies that any combination of $S$ and $T$ results in another finite order ASM in $SA_n$ with the same identity element as $S$ and $T$. Since $S_n$ is generated by the $n$-cycle $(1,2,\dots,n)$ and the transposition $(1,n)$, it follows that $\langle S, T \rangle$ is isomorphic to $S_n$.
\end{proof}

Since every finite group is isomorphic to a subgroup of a symmetric group by Cayley's Theorem, it follows that every finite group is isomorphic to a singular group of ASMs. This construction can be alternatively thought of as multiplying the following $(n+4) \times (n+4)$ idempotent matrix $\mathcal{I}_n$ on the right by any permutation matrix $P$ satisfying $P_{ii} = 1$ for $i \in \{1, 2, n+3, n+4\}$.

\[\mathcal{I}_n = \left(\begin{array}{c|c|c}
E_{22}&X_n&E_{22}\\
\hline
0&I_n&0\\
\hline
E_{11}&-X_n&E_{11}
\end{array}\right)\]

In terms of generating $S_n$ with singular ASMs of least possible size, the construction in Theorem \ref{isomorphic n+4} is best possible for certain small values of $n$. For example, this construction generates a group of $5 \times 5$ singular ASMs isomorphic to $S_1$ (the group $\{\mathcal{E}_1\}$ in Example \ref{5x5}), a group of $7 \times 7$ singular ASMs isomorphic to $S_3$ (Example \ref{2 generators}), and a group of $8 \times 8$ singular ASMs isomorphic to $S_4$ (Example \ref{8x8}). It is not optimal for $S_2$, which can be realized in $SA_5$ as in Example \ref{5x5}.

\begin{example}\label{8x8} The following singular ASMs $A$ and $B$ have order 4 and 2, and generate a group isomorphic to $S_4$ in $SA_8$.
\[\hspace{0.1cm}A = \left(\begin{array}{rr|rrrr|rr}
0 & 0 & 0 & 0 & 0 & 1 & 0 & 0 \\
0 & 1 & 0 & 0 & 0 &-1 & 0 & 1 \\
\hline
0 & 0 & 0 & 0 & 0 & 1 & 0 & 0 \\
0 & 0 & 1 & 0 & 0 & 0 & 0 & 0 \\
0 & 0 & 0 & 1 & 0 & 0 & 0 & 0 \\
0 & 0 & 0 & 0 & 1 & 0 & 0 & 0 \\
\hline
1 & 0 & 0 & 0 & 0 &-1 & 1 & 0 \\
0 & 0 & 0 & 0 & 0 & 1 & 0 & 0 
\end{array}\right)\hspace{0.1cm}
B = \left(\begin{array}{rr|rrrr|rr}
0 & 0 & 0 & 0 & 0 & 1 & 0 & 0 \\
0 & 1 & 0 & 0 & 0 &-1 & 0 & 1 \\
\hline
0 & 0 & 0 & 0 & 0 & 1 & 0 & 0 \\
0 & 0 & 0 & 1 & 0 & 0 & 0 & 0 \\
0 & 0 & 0 & 0 & 1 & 0 & 0 & 0 \\
0 & 0 & 1 & 0 & 0 & 0 & 0 & 0 \\
\hline
1 & 0 & 0 & 0 & 0 &-1 & 1 & 0 \\
0 & 0 & 0 & 0 & 0 & 1 & 0 & 0 
\end{array}\right)\]
\end{example}

\section{Groups in $SA_n$ of Lower Rank}\label{lower rank}
The idempotents constructed in Section \ref{max rank} are all diagonal extensions of the $5\times 5$ idempotent $\mathcal{E}_1$, and all have nullity 2. In contrast to this, the construction in this section involves singular groups of ASMs of unbounded size, whose identity elements are idempotents in reduced form. In particular, we show that the symmetric group $S_n$ is isomorphic to a group of ASMs of rank $2\lceil\frac{n}{2}\rceil+1$ in $SA_{4\lceil\frac{n}{2}\rceil+1}$, all in reduced form.

\begin{definition}\label{E_k}
$\mathcal{E}_k$ is defined to be the following $(4k+1)\times(4k+1)$ ASM, where each of the four $2k \times 2k$ blocks is a diagonal matrix whose diagonal entries alternate between 0 and 1.
\[\mathcal{E}_k = \left(\begin{array}{ccccccc|r|ccccccc}
 0& & & 	& & & 		& 1&		0& & &   & & & \\
 & 1& &   	& & &  		&-1&		& 1& &   & & & \\ 
 & & 0&   	& & &  		& 1&		& & 0&   & & & \\ 
 &&&\ddots 	&   &  & 		&\vdots &	& & & \ddots& & & \\ 
 & & &   	& 1& &  		&-1&		& & &   & 1& & \\ 
 & & &   	& & 0&  		& 1&		& & &   & & 0& \\ 
 & & & 	& & & 1 		&-1&		& & &   & & & 1\\
\hline
 0&0&&	\dots	& & & 		0& 1&		0& 0& &\dots & & & 0\\
\hline
1& & & 	& & & 		&-1&		1& & &  & & & \\
& 0& &  	& & & 		& 1&		& 0& &   	& & & \\ 
& & 1& 	& & & 		&-1&		& & 1&   	& & & \\ 
 &&&\ddots& & & 			& \vdots&	&&&\ddots &   &   &\\ 
& & & 	& 0& & 		& 1&		& & &   	& 0& & \\ 
& & & 	& & 1& 		&-1&		& & &   	& & 1& \\ 
& & & 	& & & 0		& 1&		& & &  	& & & 0\\
\end{array}\right)\]
\end{definition}
Note that all odd rows of $\mathcal{E}_k$ are equal, and all even rows are independent of each other and of the odd rows. Hence $\mathcal{E}_k$ is singular with rank $2k+1$. In the case $k=1$, Definition \ref{E_k} coincides with the designation of $\mathcal{E}_1$ in Example \ref{5x5}.

\begin{lemma}\label{E_k idempotent}
For all $k \in \mathbb{N}$, the matrix $\mathcal{E}_k$ is idempotent. \end{lemma}

\begin{proof}
In order to prove this lemma, we introduce the notion of a $T$-block, adapted from \cite{asbg1, asbg2}. A T-block is an $n\times n$ matrix whose non-zero entries form a (not necessarily contiguous) copy of
\[\pm\left(\begin{array}{rr}
1 & -1 \\
-1 & 1
\end{array}\right).\]

We denote by $T_{i_1,j_1;i_2,j_2}$ the $T$-block with $1$ in positions $(i_1,j_1)$ and $(i_2,j_2)$, and $-1$ in positions $(i_1, j_2)$ and $(i_2, j_1)$, where $i_1 < i_2$ and $j_1< j_2$. We note that $T_{i_1,j_1;i_2,j_2}$ is idempotent if and only if exactly one its positive entries is on the main diagonal, and that $T_{i_1,j_1;i_2,j_2}T_{i_1',j_1';i_2',j_2'}=0$ if $\{j_1,j_2\}\cap\{i_1',i_2'\}=\emptyset$. 

Let $n = 4k+1$ and $m = 2k+1$. Then $\mathcal{E}_k$ can be decomposed into $I_n$ and $2k$ T-blocks as follows.
\begin{align*}
\mathcal{E}_k = &I_n - \big(T_{1,1;m+1,m} + T_{3,3,m+3,m} + \dots + T_{m-1,m-1;n-1,m}\big)\\
&- \big(T_{2,m;m+2,m+2} + \dots + T_{m-3,n-2;n-3,n-3} + T_{m-1, n; n-1, n-1} \big)\end{align*}

Let $T$ be the sum of these T-blocks, so $\mathcal{E}_k = I_n - T$. Then $T$ is the sum of $2k$ idempotent $T$-blocks whose pairwise products are all 0, so $T^2=T$.

It follows that $(\mathcal{E}_k)^2 = (I_n - T)^2 = (I_n)^2 - 2T + T^2 = I_n - T = \mathcal{E}_k$. Therefore $\mathcal{E}_k$ is idempotent.
\end{proof}

Looking at the rows and columns of the non-zero diagonal entries, we observe that $\mathcal{E}_k$ is in reduced form, in contrast to all of our earlier examples of singular idempotent ASMs, which have nullity 2 and are diagonally extended from the $5\times 5$ reduced form $\mathcal{E}_1$ (or its transpose). Since the $(4k+1)\times(4k+1)$ matrix $\mathcal{E}_k$ has rank $2k+1$, its nullity is $2k$.

We write $\Theta_{2k}$ for the group of $(4k+1)\times (4k+1)$ permutation matrices with $1$ in the $(i,i)$ position for all odd $i$. Then $\Theta_{2k}$ is isomorphic to $S_{2k}$.

\begin{lemma}\label{PE ASM}
Let $P\in\Theta_{2k}$. Then $P\mathcal{E}_k$ is a singular ASM.
\end{lemma}
\begin{proof}
The rows of $A$ are the same as $\mathcal{E}_k$, but some of the even rows have been permuted. Therefore each row of $A$ satisfies the conditions of an ASM. Every column of $\mathcal{E}_k$ besides the central column contains one $+1$ entry as its only non-zero entry. Therefore every non-central column of $A$ satisfies the conditions of an ASM. Every even entry of the central column of $\mathcal{E}_k$ is $-1$. Therefore the central column of $A$ is equal to the central column of $\mathcal{E}_k$ and satisfies the conditions of an ASM. Since the first and last rows of $A$ are equal, $A$ is a singular ASM.
\end{proof}

\begin{lemma}\label{PQE}
Let $P,Q\in\Theta_{2k}$. Then $(P\mathcal{E}_k)(Q\mathcal{E}_k)=PQ\mathcal{E}_k.$
\end{lemma}
\begin{proof}
As in the proof of Lemma \ref{E_k idempotent}, $\mathcal{E}_k = I - T$, where $T$ consists of $2k$ T-blocks occupying different rows and columns of $T$, and each T-block has a $+1$ in an odd row/column along the diagonal. The odd row also contains a $-1$ entry in the central column, and the odd column also contains a $-1$ entry in an even row. Therefore the final non-zero entry of the T-block is a $+1$ in the central column of this even row.

So $AB = (P\mathcal{E}_k)(Q\mathcal{E}_k) = P(I - T)Q(I - T) = (P-PT)(Q-QT) = PQ - PQT - PTQ + PTQT$.

Since multiplying a matrix on the right by $Q$ permutes only even columns of that matrix, and all even columns of $T$ are empty, it follows that $TQ = T$. Therefore $AB = PQ - PQT - PT + PT^2$.

As in the proof of Lemma \ref{E_k idempotent}, $T^2 = T$. Therefore 
\[AB = PQ - PQT - PT = PT = PQ(I-T) = (PQ)\mathcal{E}_k.\]
\end{proof}

\begin{theorem}\label{PEk}
Let $P$ be a $(4k+1)\times(4k+1)$ permutation matrix of multiplicative order $n$ such that $P_{ii} = 1$ for all odd $i$. Then $A = P\mathcal{E}_k$ has order $n$ and $A^m$ is a singular ASM for all $m \in \{1,2,\dots,n\}$.
\end{theorem}

\begin{proof}
By Lemma \ref{PE ASM}, $A^m$ is a singular ASM for all $m$. By Lemma \ref{PQE}, $A = P^m\mathcal{E}_k$ for all $m$. Since $\mathcal{E}_k$ is idempotent and $P$ has order $n$, it follows that $A$ has order $n$.
\end{proof}

It follows from Lemmas \ref{PE ASM} and \ref{PQE} and Theorem \ref{PEk} that the mapping from $\Theta_{2k}$ into $SA_{4k+1}$ is a homomorphism of groups, that preserves element orders and is therefore injective. Thus $\Theta_{2k}\mathcal{E}_k$ is a copy of the symmetric group $S_{2k}$, consisting of matrices of rank $2k$ in $SA_{4k+1}$.

\begin{example}
The following $9\times9$ singular ASMs $A$ and $B$ have order 3 and 4, respectively.
\[A = \left(\begin{array}{cccc r cccc}
0&0&0&0&	1	&0&0&0&0\\
1&0&0&0&	-1	&1&0&0&0\\
0&0&0&0&	1	&0&0&0&0\\
0&1&0&0&	-1	&0&1&0&0\\

0&0&0&0&	1	&0&0&0&0\\

0&0&0&1&	-1	&0&0&0&1\\
0&0&0&0&	1	&0&0&0&0\\
0&0&1&0&	-1	&0&0&1&0\\
0&0&0&0&	1	&0&0&0&0\\
\end{array}\right)
\hspace{1cm}
B = \left(\begin{array}{cccc r cccc}
0&0&0&0&	1	&0&0&0&0\\
0&0&1&0&	-1	&0&0&1&0\\
0&0&0&0&	1	&0&0&0&0\\
0&1&0&0&	-1	&0&1&0&0\\

0&0&0&0&	1	&0&0&0&0\\

0&0&0&1&	-1	&0&0&0&1\\
0&0&0&0&	1	&0&0&0&0\\
1&0&0&0&	-1	&1&0&0&0\\
0&0&0&0&	1	&0&0&0&0\\
\end{array}\right)\]
The group $\langle A, B \rangle$ is isomorphic to $S_4$, and each element is a singular ASM.
\end{example}

\begin{theorem}\label{group representation}
$S_n$ is isomorphic to a group in $SA_{4\lceil\frac{n}{2}\rceil+1}$ consisting of matrices of rank $2\lceil\frac{n}{2}\rceil+1$.
\end{theorem}

\begin{proof}
$S_n$ is generated by the $n$-cycle $(1,2,\dots,n)$ and the transposition $(1,n)$. Let $k$ by the least integer such that $2k \geq n$. So $k = \lceil\frac{n}{2}\rceil$. Let $P_n$ be the $(4k+1)\times(4k+1)$ permutation matrix resulting from cycling the first $n$ even rows of the $(4k+1)\times(4k+1)$ identity matrix, and let $P_2$ be the permutation matrix resulting from swapping the first 2 even rows.

Then, by Lemma \ref{PE ASM} and Theorem \ref{PEk}, $S = P_n\mathcal{E}_k$ and $T = P_2\mathcal{E}_k$ generate a copy of 
 $S_n$ in $SA_{4k+1}$.
\end{proof}

\begin{example}
Since $4\lceil\frac{3}{2}\rceil+1 = 9$, a group in $SA_{9}$ isomorphic to $S_3$ is generated by the following matrices.
\[S = \left(\begin{array}{cccc r cccc}
0&0&0&0&	1	&0&0&0&0\\
1&0&0&0&	-1	&1&0&0&0\\
0&0&0&0&	1	&0&0&0&0\\
0&1&0&0&	-1	&0&1&0&0\\

0&0&0&0&	1	&0&0&0&0\\

0&0&0&1&	-1	&0&0&0&1\\
0&0&0&0&	1	&0&0&0&0\\
0&0&1&0&	-1	&0&0&1&0\\
0&0&0&0&	1	&0&0&0&0\\
\end{array}\right)
\hspace{1cm}
T = \left(\begin{array}{cccc r cccc}
0&0&0&0&	1	&0&0&0&0\\
0&0&0&1&	-1	&0&0&0&1\\
0&0&0&0&	1	&0&0&0&0\\
0&1&0&0&	-1	&0&1&0&0\\

0&0&0&0&	1	&0&0&0&0\\

1&0&0&0&	-1	&1&0&0&0\\
0&0&0&0&	1	&0&0&0&0\\
0&0&1&0&	-1	&0&0&1&0\\
0&0&0&0&	1	&0&0&0&0\\
\end{array}\right)\]
\end{example}

Note that the constructions outlined in this paper can be combined to generate more examples of singular ASMs of finite order. 

\begin{example}The following $9 \times 9$ singular ASM $A$ has order 4 and arises from the construction of Theorem \ref{PEk}. The following $5\times5$ permutation matrix $P$ has order 5.
\[A = \left(\begin{array}{cccc r cccc}
0&0&0&0&	1	&0&0&0&0\\
0&0&1&0&	-1	&0&0&1&0\\
0&0&0&0&	1	&0&0&0&0\\
0&1&0&0&	-1	&0&1&0&0\\

0&0&0&0&	1	&0&0&0&0\\

0&0&0&1&	-1	&0&0&0&1\\
0&0&0&0&	1	&0&0&0&0\\
1&0&0&0&	-1	&1&0&0&0\\
0&0&0&0&	1	&0&0&0&0\\
\end{array}\right)
\hspace{1cm}
P = \left(\begin{array}{ccccc}
0&0&0&0&1\\
1&0&0&0&0\\
0&1&0&0&0\\
0&0&1&0&0\\
0&0&0&1&0\\
\end{array}\right)\]

Combining $A$ with $P$ using the construction of Theorem \ref{block perm orders} results in $B \in SA_{13}$ of order 20.
\[B = \left(\begin{array}{cccc |r| cccc}
0&0&0&0&	1	&0&0&0&0\\
0&0&1&0&	-1	&0&0&1&0\\
0&0&0&0&	1	&0&0&0&0\\
0&1&0&0&	-1	&0&1&0&0\\
\hline
0&0&0&0&	P	&0&0&0&0\\
\hline
0&0&0&1&	-1	&0&0&0&1\\
0&0&0&0&	1	&0&0&0&0\\
1&0&0&0&	-1	&1&0&0&0\\
0&0&0&0&	1	&0&0&0&0\\
\end{array}\right)
= \left(\begin{array}{cccc |ccccr| cccc}
0&0&0&0&	0&0&0&0&1	&0&0&0&0\\
0&0&1&0&	0&0&0&0&-1	&0&0&1&0\\
0&0&0&0&	0&0&0&0&1	&0&0&0&0\\
0&1&0&0&	0&0&0&0&-1	&0&1&0&0\\
\hline
0&0&0&0&	0&0&0&0&1	&0&0&0&0\\
0&0&0&0&	1&0&0&0&0	&0&0&0&0\\
0&0&0&0&	0&1&0&0&0	&0&0&0&0\\
0&0&0&0&	0&0&1&0&0	&0&0&0&0\\
0&0&0&0&	0&0&0&1&0	&0&0&0&0\\

\hline
0&0&0&1&	0&0&0&0&-1	&0&0&0&1\\
0&0&0&0&	0&0&0&0&1	&0&0&0&0\\
1&0&0&0&	0&0&0&0&-1	&1&0&0&0\\
0&0&0&0&	0&0&0&0&1	&0&0&0&0\\
\end{array}\right)\]
\end{example}

\section{Kronecker Product of Groups of ASMs}
Let $A$ and $B$ be matrices of size $m$ and $n$, respectively. The \emph{Kronecker product} $A \otimes B$ is the $(mn) \times (mn)$ matrix consisting of $n \times n$ blocks arranged in an $m \times m$ array such that the $(i,j)$ block is equal to $a_{ij}B$. For background on the Kronecker product, we refer to Section 4.3 of \cite{kronecker}.

For example, recall $\mathcal{E}_1 = \left(\begin{array}{rrrrr}
0 & 0 & 1 & 0 & 0 \\
0 & 1 & -1 & 0 & 1 \\
0 & 0 & 1 & 0 & 0 \\
1 & 0 & -1 & 1 & 0 \\
0 & 0 & 1 & 0 & 0
\end{array}\right)$. Then $\mathcal{E}_1 \otimes \mathcal{E}_1 = \left(\begin{array}{c|c|c|c|c}
0&0&\mathcal{E}_1&0&0\\
\hline
0&\mathcal{E}_1&-\mathcal{E}_1&0&\mathcal{E}_1\\
\hline
0&0&\mathcal{E}_1&0&0\\
\hline
\mathcal{E}_1&0&-\mathcal{E}_1&\mathcal{E}_1&0\\
\hline
0&0&\mathcal{E}_1&0&0\\
\end{array}\right)$.

Note that $\rank(A\otimes B) = \rank(A)\rank(B)$ \cite{kronecker}. So $\mathcal{E}_1 \otimes \mathcal{E}_1$ has size 25 and rank 9.

\begin{lemma}
    If $A$ and $B$ are ASMs, then $A \otimes B$ is an ASM.
\end{lemma}
\begin{proof}
    Consider any row of $A \otimes B$. This row intersects at least one non-zero block, the first of which is itself an ASM, so these non-zero entries alternate in sign, beginning and ending with $+1$. If the row intersects a second non-zero block, it is the negative of an ASM, and so these non-zero entries alternate in sign, beginning and ending with $-1$. The non-zero entries continue to alternate like this, and the last non-zero block that this row intersects is an ASM, so the non-zero entries in the whole row alternate in sign, beginning and ending with $+1$. The same is true for any column, and therefore $A \otimes B$ is an ASM.
\end{proof}

Note that $(A \otimes B)(C \otimes D) = (AC)\otimes(BD)$ for all matrices $A, B, C, D$ for which $AC$ and $BD$ are defined \cite{kronecker}. This means that if $A$ and $B$ have multiplicative orders $a$ and $b$, the order of $A \otimes B$ is $\lcm(a,b)$. So, for example, if $A$ and $B$ are the singular $7 \times 7$ ASMs of order 3 and 2 in Example \ref{2 generators}, then $A \otimes B$ is a $49 \times 49$ singular ASM of order 6.

\[\left(\begin{array}{rrrrrrr}
0 & 0 & 0 & 0 & 1 & 0 & 0 \\
0 & 1 & 0 & 0 &-1 & 0 & 1 \\
0 & 0 & 0 & 0 & 1 & 0 & 0 \\
0 & 0 & 1 & 0 & 0 & 0 & 0 \\
0 & 0 & 0 & 1 & 0 & 0 & 0 \\
1 & 0 & 0 & 0 &-1 & 1 & 0 \\
0 & 0 & 0 & 0 & 1 & 0 & 0 
\end{array}\right) 
\otimes \left(\begin{array}{rrrrrrr}
0 & 0 & 0 & 0 & 1 & 0 & 0 \\
0 & 1 & 0 & 0 &-1 & 0 & 1 \\
0 & 0 & 0 & 0 & 1 & 0 & 0 \\
0 & 0 & 0 & 1 & 0 & 0 & 0 \\
0 & 0 & 1 & 0 & 0 & 0 & 0 \\
1 & 0 & 0 & 0 &-1 & 1 & 0 \\
0 & 0 & 0 & 0 & 1 & 0 & 0 
\end{array}\right)
= \left(\begin{array}{r|r|r|r|r|r|r}
0 & 0 & 0 & 0 & B & 0 & 0 \\
\hline
0 & B & 0 & 0 &-B & 0 & B \\
\hline
0 & 0 & 0 & 0 & B & 0 & 0 \\
\hline
0 & 0 & B & 0 & 0 & 0 & 0 \\
\hline
0 & 0 & 0 & B & 0 & 0 & 0 \\
\hline
B & 0 & 0 & 0 &-B & B & 0 \\
\hline
0 & 0 & 0 & 0 & B & 0 & 0  
\end{array}\right)\]

If $G$ and $H$ are groups in $SA_m$ and $SA_n$ respectively, and we define $G \otimes H = \{g \otimes h \::\: g \in G, h \in H\}$, then $G \otimes H$ is a group in $SA_{mn}$. Multiplication in $G \otimes H$ satisfies $(g_1 \otimes h_1)(g_2 \otimes h_2) = g_1g_2 \otimes h_1h_2$, and therefore $G \otimes H$ is isomorphic to the direct product $G \times H$.
Thus the Kronecker product provides a mechanism for realizing direct products as groups of singular ASMs.


\section{Conclusions}
In this final section, we suggest some directions for further research on singular groups of ASMs.

The multiplicative order of every singular ASM presented in this paper coincides with that of some permutation matrix of the same size. In \cite{finiteASM}, it is shown that there exist nonsingular ASMs that generate finite cyclic groups whose order is not equal to that of any permutation in the corresponding symmetric group. We pose the following problem. 

\begin{problem}
    Does there exist a cyclic group in $SA_n$ whose order is not that of any cyclic subgroup of $S_n$?
\end{problem}

Example \ref{8x8} uses the construction of Theorem \ref{isomorphic n+4} to generate a group isomorphic to $S_4$ using matrices in $SA_8$. By exhaustive search, we know that this construction is optimal in the sense that $S_4$ cannot be realised as a group of singular ASMs of smaller size. This construction is also optimal for $S_3$, but not for $S_2$ which is a cyclic group. We conjecture that this construction is optimal for all $n > 2$.

\begin{conjecture}
For $n>2$, the the least $k$ for which $SA_k$ contains a group isomorphic to $S_n$ is $k=n+4$.
\end{conjecture}

Each example and construction outlined for an idempotent in $SA_n$ has had all of its negative entries in a single column (or row). We have developed methods for constructing wider classes of idempotent ASMs, which will be the subject of a forthcoming paper.

\end{document}